\documentclass{amsart}

\theoremstyle{plain}
\newtheorem{Theorem}{Theorem}
\newtheorem{Corollary}[Theorem]{Corollary}
\newtheorem{Proposition}[Theorem]{Proposition}

\theoremstyle{definition}

\theoremstyle{remark}
\newtheorem*{Remark}{Remark}
\newtheorem*{Notation}{Notation}

\DeclareMathOperator{\Bs}{Bs}
\DeclareMathOperator{\Div}{Div}
\DeclareMathOperator{\Pic}{Pic}
\DeclareMathOperator{\Supp}{Supp}
\DeclareMathOperator{\mult}{mult}
\DeclareMathOperator{\ord}{ord}

\newcommand{\Z}{\ensuremath{\mathbb{Z}}}
\newcommand{\Q}{\ensuremath{\mathbb{Q}}}

\newcommand{\C}{\ensuremath{\mathbb{C}}}

\newcommand{\Prsp}{\ensuremath{\mathbb{P}}}

\newcommand{\Ocal}{\mathcal{O}}
\newcommand{\Ecal}{\mathcal{E}}
\newcommand{\I}{\mathcal{I}}

\newcommand{\MM}{\frak{m}}

\newcommand{\union}{\cup}
\newcommand{\inters}{\cap}
\newcommand{\tens}{\otimes}
\newcommand{\lra}{\longrightarrow}

\newcommand{\qle}{\underset{\Q}{\sim}}  
\newcommand{\eqdef}{\overset{\text{def}}{=}}

\newcommand{\vv}{\vec{v}}

\newcommand{\ulc}{\lceil}
\newcommand{\urc}{\rceil}
\newcommand{\llc}{\lfloor}
\newcommand{\lrc}{\rfloor}

\newcommand{\rup}[1]{\ulc #1 \urc}	
\newcommand{\rdn}[1]{\llc #1 \lrc}	
\newcommand{\frp}[1]{\{ #1 \}}		

\begin{document}



\title[Very ampleness on surfaces with boundary] {Very ampleness of 
		adjoint linear systems on smooth surfaces with boundary}

\author{Vladimir Ma\c{s}ek}
\address{Department of Mathematics, Box 1146, Washington University,
		St. Louis, MO 63130}
\email{vmasek@math.wustl.edu}


\begin{abstract}
  Let $M$ be a \Q-divisor on a smooth surface over $\C\,$.
  In this paper we give criteria for very ampleness of the adjoint of
  $\rup{M}$, the round-up of $M$. (Similar results for global generation
  were given by Ein and Lazarsfeld and used in their proof of Fujita's 
  Conjecture in dimension 3.) In \S 4 we discuss an example which suggests 
  that this kind of criteria might also be useful in the study of linear
  systems on surfaces.
\end{abstract}

\maketitle




\setcounter{section}{-1}
\subsection*{Contents}

\begin{enumerate}
  \item[0.] Introduction 
  \item[1.] Base-point-freeness
  \item[2.] Separation of points
  \item[3.] Separation of tangent directions
  \item[4.] Example
\end{enumerate}

\subsection*{Notations}

\begin{tabbing}
99\=9999999999\=9999999999999999999999999999\kill
  \>$\ulc \cdot \urc$ \> round-up               \\
  \>$\llc \cdot \lrc$ \> round-down             \\
  \>$\{ \cdot \}$     \> fractional part        \\
  \>$f^{-1}D$         \> strict transform (proper transform) \\
  \>$f^*D$            \> pull-back (total inverse image) \\
  \>$PLC$	      \> partially log-canonical (Definition 1.7) \\
  \>$\,\equiv$        \> numerical equivalence  \\
  \>$\,\sim$          \> linear equivalence     \\
  \>$\,\qle$          \> \Q-linear equivalence  \\
\end{tabbing}


\section{Introduction}

Let $S$ be a nonsingular projective surface over $\C\,$, and let $H$ be a
given line bundle on $S$. Consider the following natural questions regarding
the complete linear system $|H|$:
\emph{ 
  \begin{enumerate}
      \item[(1)] Compute $\dim |H|$.
      \item[(2)] Is $|H|$ base-point-free?
      \item[(3)] Is $|H|$ very ample?
  \end{enumerate}
}
The answer to (1) is usually given in two parts: the Riemann-Roch theorem
computes $\chi(S,H)$, and then we need estimates for $h^i(S,H),\; i>0$. In
particular, we may ask the following question related to (1):
\emph{
  \begin{enumerate}
    \item[($1'$)] When are $h^1(S,H) \text{ and } h^2(S,H)$ equal to zero?
  \end{enumerate}
}
One classical answer to ($1'$) is provided by Kodaira's vanishing theorem:
if $L$ is any ample line bundle on $S$, then $h^i(S,-L)=0$ for all $i<2$;
therefore, by Serre duality, we have $h^i(S, K_S+L)=0$ for all $i > 0$. To
answer ($1'$), write $H=K_S+L$ (thus defining $L$ as $H-K_S$); if $L$ is
ample, then $h^i(S,H)=0$ for all $i>0$.

\vspace{6pt}

For questions (2) and (3), Reider \cite{rei} gave an answer which again
considers $H$ in the form of an adjoint line bundle, $H=K_S+L$:

{\bf Proposition} (cf. \cite[Theorem 1]{rei}){\bf .}  \emph{
If $L$ is a line bundle on $S$, 
$L^2 \geq 5$ and $L \cdot C \geq 2$ for every curve $C \subset S$, then 
$|K_S+L|$ is base-point-free. If $L^2 \geq 10$ and $L \cdot C \geq 3$ for
every curve $C$, then $|K_S+L|$ is very ample.
}

\vspace{6pt}

We note here that Kodaira's theorem holds in all dimensions. Reider's
criterion was tentatively extended in higher dimensions in the form of
Fujita's conjecture (\cite{fuj}): if $X$ is a smooth projective variety
of dimension $n$, and $L$ is an ample line bundle on $X$, then $|K_X+mL|$ is
base-point-free for $m\geq n+1$ and very ample for $m\geq n+2$. Fujita's
conjecture for base-point-freeness was proved in dimension 3 by Ein and
Lazarsfeld (\cite{el}) and in dimension 4 by Kawamata (\cite{kaw}); more
precise statements, which resemble Reider's criterion more closely, were also
obtained. Very ampleness, however, is still open, even in dimension 3.

Kodaira's vanishing theorem and Reider's criterion are already very useful as
stated. However, the applicability of Kodaira's theorem was greatly extended,
first on surfaces, by Mumford, Ramanujam, Miyaoka, and then in all dimensions
by Kawamata and Viehweg, as follows. First, the ampleness condition for $L$ can
be relaxed to $L\cdot C \geq 0$ for every curve $C$ and $L^2>0$ ($L$ \emph{nef}
and \emph{big}). Second, and most important, assume that $L$ itself is not nef
and big, but there is a nef and big \Q-divisor $M$ on $S$ ($M \in
\Div(S)\tens\Q$) such that $L = \rup{M}$ (i.e. $L-M$ is an effective \Q-divisor
$B$ whose coefficients are all $<1$). Then we have $h^i(S,K_S+L)=0$ for all
$i>0$, just as in Kodaira's theorem. 
(\Q-divisors were first considered in this context in connection with the
Zariski decomposition of effective divisors.)

In dimension $\geq 3$, the Kawamata--Viehweg vanishing theorem requires an
extra hypothesis (the irreducible components of $\Supp(B)$ must cross
normally); however, Sakai remarked that for surfaces this extra hypothesis
is not necessary (see Proposition 1.2.1 in \S 1).

\vspace{6pt}

For base-point-freeness (question (2) above), Ein and Lazarsfeld (\cite{el})
proved a similar extension of Reider's criterion, expressing $H$ as
$K_S+\rup{M}$ for a \Q-divisor $M$ on $S$; if $M^2>4$ and $M\cdot C \geq 2$
for every curve $C$, then $|H|$ is base-point-free. They used this result in
their proof of Fujita's conjecture for base-point-freeness in dimension 3.
(In fact they used a more precise local version, involving the local
multiplicities of $B=L-M$; see \S 1 below).

In this paper we give criteria for very ampleness of linear systems of the
form $|K_S+B+M|$, $B=\rup{M}-M$, as above. In particular, we prove the
following result:

\begin{Theorem}
  Let $S$, $B$ and $M$ be as above, and assume that 
    \begin{enumerate}
      \item[(0.1)] $M^2 > 2 (\beta_2)^2$,
      \item[(0.2)] $M \cdot C \geq 2\beta_1$ for every irreducible curve
      		   $C \subset S$, 				\newline
      		   where $\beta_2$, $\beta_1$ are positive
      		   numbers satisfying the following inequalities:
      \item[(0.3)] $\beta_2 \geq 2$,
      \item[(0.4)] $\beta_1 \geq \dfrac{\beta_2}{\beta_2-1}$.
    \end{enumerate}
  Then $|K_S+B+M|$ is very ample.
\end{Theorem}

An immediate consequence of Theorem 1 is the following:
\begin{Corollary}
  Assume that (S,B) is as before, and $M$ is an ample \Q-divisor on $S$ such
  that $B=\rup{M}-M$, $M^2 > (2+\sqrt{2})^2$, and $M \cdot C > 2+\sqrt{2}$ for
  every curve $C \subset S$. Then $|K_S+B+M|$ is very ample. In particular, 
  if $A$ is an ample divisor (with integer coefficients) on $S$, then
  $|K_S + \rup{aA}|$ is very ample for every $a \in \Q$, $a > 2+\sqrt{2}$.
\end{Corollary}
Note that Reider's criterion implies only that $|K_S+aA|$ is very ample for
every \emph{integer} $a \geq 4$.

\vspace{6pt}

As in \cite[\S 2]{el} (where the analogue for base-point-freeness was
proved), we prove a local version of Theorem 1, with the numerical
conditions on $M$ relaxed in terms of local multiplicities of $B$.

\vspace{6pt}

As we mentioned earlier, the result for base-point-freeness on surfaces with
boundary (i.e. for \Q-divisors $M$) was used in \cite{el} in the proof of
Fujita's Conjecture in dimension 3. Similarly, we expect that the proof of
the analogous result for very ampleness in dimension 3 will use very ampleness
for \Q-divisors on surfaces. However, a natural and interesting question is
whether or not the results for \Q-divisors on surfaces have any useful
applications to the study of linear systems on surfaces. An example we discuss
in \S 4 seems to indicate an affirmative answer. While the results proved in
\S 4 can be obtained with other methods, our example shows how our \Q-Reider
theorem extends the applicability of Reider's original result in the same way
the Kawamata--Viehweg vanishing theorem extends the range of applicability of
Kodaira's vanishing theorem. The usefulness of considering local
multiplicities of $B$ is also evident in this example.

\vspace{6pt}

The paper is divided as follows: \S 1 is devoted to base-point-freeness.
The results discussed in this section, with one exception, were proved in
\cite{el}; I include a (slightly modified) proof to fix the ideas and notations
for the later sections. As one might expect, separation of points is relatively
easy (at least in principle); it is discussed in \S 2. Then we move on to
separation of tangent directions in \S 3. This part is surprisingly delicate;
in particular the ``multiplier ideal'' method of Ein--Lazarsfeld, or 
Kawamata's  equivalent ``log-canonical threshold'' formalism, do not work in
this context. We explain the geometric contents of our method in the
beginning of \S 3. Theorem 1 follows from Proposition 4 in \S 2 and
Proposition 5 in \S 3. Finally, \S 4 contains the example mentioned earlier.

\vspace{6pt}

The author is grateful to L. Ein, R. Lazarsfeld, S. Lee, and N. Mohan Kumar
for their many useful suggestions.


\section{Base-point-freeness}

{\bf (1.1)} Let $S$ be a smooth projective surface over \C, and 
$B=\sum b_i C_i$ a fixed effective \Q-divisor on $S$ with $0 \leq b_i < 1$
for all $i$. (The pair $(S,B)$ is sometimes called a ``surface with
boundary'', whence the title of this paper.) Let $M$ be a \Q-divisor on $S$
such that $B+M$ has integer coefficients.

{\bf We assume throughout this paper that $M$ is nef and big,} i.e. that
$M \cdot C \geq 0$ for every curve $C \subset S$ and $M^2 > 0$.

\vspace{6pt}

{\bf (1.2)} For convenience, we gather here two technical results which we use
time and again in our proofs.

{\bf (1.2.1)} We use the following variants of the Kawamata--Viehweg vanishing
theorem, which hold on smooth surfaces:

{\bf Theorem.} {\bf (a)} (cf. \cite[Lemma 1.1]{el}) Let $S$ be a smooth
projective surface over \C, and let $M$ be a nef and big \Q-divisor on $S$.
Then 
\begin{equation*}
  H^i(S, K_S+\rup{M})=0, \qquad \forall i > 0.
\end{equation*}

{\bf (b)} (cf. \cite[Lemma 2.4]{el}) Assume moreover that $C_1, \dots, C_k$
are distinct irreducible curves on $S$ which have integer coefficients in $M$.
Assume that $M \cdot C_j > 0$ for all $j = 1, \dots, k$. Then 
\begin{equation*}
  H^i(S, K_S+\rup{M}+C_1+\cdots+C_k)=0, \qquad \forall i > 0.
\end{equation*}							\qed

\vspace{6pt}

{\bf (1.2.2)} We use the following criterion for base-point-freeness,
respectively very ampleness, on a complete Gorenstein curve (cf. \cite{har}):

{\bf Proposition.} Let $D$ be a Cartier divisor on the integral
projective Gorenstein curve $C$. Then:

{\bf (a)} $\deg(D) \geq 2 \implies$ the complete linear system $|K_C+D|$
	  is base-point-free;

{\bf (b)} $\deg(D) \geq 3 \implies |K_C+D|$ is very ample.

\begin{proof}
  See \cite[\S 1]{har} for the relevant definitions (generalized divisors on
  $C$, including $0$-dimensionals subschemes; degree; etc.)

  We prove (b); the proof of (a) is similar. By \cite[Proposition 1.5]{har},
  it suffices to show that $h^0(C, K_C+D-Z)=h^0(C, K_C+D)-2$ for every
  $0$-dimensional subscheme $Z \subset C$ of length $2$. Consider the exact
  sequence: 
 \begin{equation*}
    0 \lra \Ocal_C(K_C+D-Z) \lra \Ocal_C(K_C+D) \lra
		\Ocal_C(K_C+D) \tens \Ocal_Z \lra 0.
 \end{equation*}
  As $\Ocal_C(K_C+D)\tens\Ocal_Z \cong \Ocal_Z$ has length $2$, the conclusion
  will follow from the vanishing of $H^1(C, K_C+D-Z)$. By Serre duality (cf.
  \cite[Theorem 1.4]{har}), $H^1(C, K_C+D-Z) \cong H^0(C, Z-D)$, and
  $H^0(C, Z-D)=0$ due to $\deg(Z-D) = 2 - \deg(D) < 0$.
\end{proof}

\vspace{6pt}

{\bf (1.3)} Fix a point $p \in S$. In this section we give sufficient
conditions for $|K_S+B+M|$ to be free at $p$.

\vspace{3pt}

\noindent {\bf (1.3.1)} {\it Notation.}
  $\quad \mu = \ord_p(B) \eqdef \sum b_i \cdot \mult_p(C_i) \qquad
  	(B = \sum b_i C_i).$

\vspace{4pt}

\begin{Proposition}
  $|K_S+B+M|$ is free at $p$ in each of the following cases:
  \begin{enumerate}
    \item[{\bf 1.}] $\mu \geq 2$;
    \item[{\bf 2.}] $0 \leq \mu < 2$; $M^2 > (\beta_2)^2$, $M \cdot C \geq
              \beta_1$ for every irreducible curve $C \subset S$ such that
              $p \in C$, where $\beta_2$, $\beta_1$ are positive numbers which
              satisfy the inequalities:
  \end{enumerate} 
  \begin{align*}
    \beta_2 &\geq 2-\mu ,    			\tag{1.3.2}		\\
    \beta_1 &\geq \min \left\{ (2-\mu) ; \frac{\beta_2}{\beta_2-(1-\mu)} 
    			\right\} .		\tag{1.3.3}
  \end{align*}
\end{Proposition}

\begin{Remark}
  Explicitly, the minimum in (1.3.3) is given by: 
    \begin{equation*}
      \min \left\{ (2-\mu) ; \frac{\beta_2}{\beta_2-(1-\mu)} \right\} =
          \begin{cases}
              2-\mu 			     & \text{if $1 \leq \mu < 2$} \\
              \dfrac{\beta_2}{\beta_2-(1-\mu)} & \text{if $0 \leq \mu < 1$.}
          \end{cases}
    \end{equation*}
  In other words, when $0 \leq \mu < 2$, the inequalities $\beta_2 \geq 2-\mu$
  and $\beta_1\geq 2-\mu$ suffice. When $\mu < 1$ the inequality for $\beta_1$ 
  can be relaxed to 
    \begin{equation*}
      \beta_1 \geq \frac{\beta_2}{\beta_2-(1-\mu)} ;  \tag{1.3.4}
    \end{equation*}
  this last part (which is useful in applications, cf. \S 4) is not
  contained in \cite{el}.
\end{Remark}

\vspace{6pt}

\noindent {\bf Proof of Proposition 3.}

\vspace{4pt}

{\bf (1.4)} Let $f : S_1 \to S$ be the blowing-up of $S$ at $p$, and let 
$E \subset S_1$ be the exceptional divisor of $f$. We have
$f^*B=f^{-1}B+\mu E$; $\rdn{f^{-1}B}=0$, and therefore 
\begin{equation*}
  \begin{split}
    K_{S_1}+\rup{f^*M} &= f^*K_S+E+\rup{f^*(B+M)-f^*B}  	\\
  		       &= f^*(K_S+B+M)+E-\rdn{f^*B}		\\   
  		       &= f^*(K_S+B+M)-(\rdn{\mu}-1)E.    
  \end{split}					\tag{1.4.1} 
\end{equation*}

\vspace{4pt}

{\bf (1.5)} If $\mu \geq 2$, then $p \notin \Bs |K_S+B+M|$. Indeed, in this
case $t=\rdn{\mu}-1$ is a positive integer; since $f^*M$ is nef and big on
$S_1$, the vanishing theorem (1.2.1)(a) yields 
\begin{equation*}
  H^1(S_1, K_{S_1}+\rup{f^*M})=0,		\tag{1.5.1}
\end{equation*}
and therefore (using (1.4.1) and the projection formula) 
\begin{equation*}
  H^1(S, \Ocal_S(K_S+B+M) \tens {\MM_p}^t)=0,	\tag{1.5.2}
\end{equation*}
where $\MM_p$ is the maximal ideal of $\Ocal_S$ at $p$. The conclusion 
follows from the surjectivity of the restriction map 
\begin{equation*}
  H^0(S, K_S+B+M) \lra H^0(S, \Ocal_S(K_S+B+M)\tens
  	\Ocal_S/{\MM_p}^t) \cong \Ocal_S/{\MM_p}^t .
\end{equation*}

\vspace{4pt}

{\bf (1.5.3)} \emph{Remark:} In fact we proved that $|K_S+B+M|$ separates
$s$-jets at $p$, if $\mu \eqdef \ord_p(B) \geq s+2$.

\vspace{9pt}

{\bf (1.6)} Now assume that $\mu < 2$, and $M^2 > (\beta_2)^2$ with
$\beta_2 \geq 2-\mu$, etc.

\vspace{6pt}

{\bf (1.6.1) Claim:} We can find an effective \Q-divisor $D$ on $S$
  such that $\ord_p(D) = 2-\mu$ and $D \qle tM$ for some $t \in \Q\,$,
  $0 < t < \dfrac{2-\mu}{\beta_2}$. ($\qle$ denotes \Q-linear equivalence,
  i.e.  $mD$ and $mtM$ have integer coefficients and are linearly equivalent
  for some suitably large and divisible integer $m$.)

\vspace{4pt}

\emph{Proof of (1.6.1):} By Riemann--Roch, $\dim |nM|$ grows like
  $\frac{M^2}{2}n^2 > \frac{(\beta_2)^2}{2}n^2$ for $n$ sufficiently large
  and divisible (such that $nM$ has integer coefficients). 
  Since $\dim (\Ocal_{S,p}/{\MM_p}^n)$ grows like $\frac{n^2}{2}$,
  for suitable $n$ we can find $G \in |nM|$ with $\ord_p(G) > \beta_2 n$.

  Take $D=rG$, $r=\dfrac{2-\mu}{\ord_p(G)}$; then $\ord_p(D) = 2-\mu$, and
  $D \qle tM$ for $t=rn < \dfrac{2-\mu}{\beta_2 n}n = \dfrac{2-\mu}{\beta_2}$.
			\qed

\vspace{6pt}

Note that $\dfrac{2-\mu}{\beta_2} \leq 1$, by (1.3.2), so that $t < 1$;
therefore $M-D \qle (1-t)M$ is still nef and big.

\vspace{8pt}

{\bf (1.7)} Recall that $B = \sum b_i C_i$, for distinct irreducible curves
$C_i \subset S$. Write $D = \sum d_i C_i$ (we allow some coefficients $b_i$
and $d_i$ to be zero); $d_i \in \Q$, $d_i \geq 0$, and
	$\ord_p(D) = \sum d_i \cdot \mult_p(C_i) = 2-\mu$.

Let $D_i = f^{-1}C_i \subset S_1$ be the strict transform of $C_i$; then
  $f^*B = \sum b_i D_i + \mu E$, $f^*D = \sum d_i D_i + (2-\mu) E$, 
  $K_{S_1}=f^*K_S+E$, and 
\begin{equation*}
  K_{S_1} - f^*(K_S+B+D) = -E - \sum (b_i + d_i) D_i.
\end{equation*}

\vspace{4pt}

\emph{Definition.} $(S,B,D)$ is {\bf partially log-canonical at $p$}
  ($PLC$ at $p$) if $-(b_i + d_i) \geq -1$ (i.e.  $b_i+d_i \leq 1$) for every
  $i$ such that $p \in C_i$. (The general definition requires the coefficient
  of $E$ to be $\geq -1$, too; in our case that coefficient is equal to $-1$.)

Note that $PLC$ is not the same as \emph{log-canonical} (cf.
\cite[Definition 0-2-10]{kmm}), because $f$ is not an embedded resolution
of $(S, B+D)$.

\vspace{6pt}

{\bf (1.8)} If $(S,B,D)$ is $PLC$ at $p$, then the proof is almost as simple
as in the case $\mu \geq 2$:
\begin{equation*}
  \begin{split}
     K_{S_1}+\rup{f^*(M-D)} &= f^*K_S+E+f^*(B+M)-\rdn{f^*(B+D)}		\\
               &= f^*(K_S+B+M)+E-2E-\sum\rdn{b_i+d_i}D_i		\\
               &= f^*(K_S+B+M)-E-\textstyle\sum ' D_i -N_1,
  \end{split}						\tag{1.8.1}
\end{equation*}
  where $\sum ' D_i$ extends over those $i$ for which $p \in C_i$ and
  $b_i+d_i=1$ (if any), and $N_1$ is an effective divisor supported away
  from $E$.

$f^*M \cdot D_i = M \cdot C_i > 0$ if $p \in C_i$; therefore (1.2.1)(b) yields: 
\begin{equation*}
  H^1(S_1, f^*(K_S+B+M)-E-N_1) = 0.			\tag{1.8.2}
\end{equation*}

Arguing as in (1.5), we can show that $p \notin \Bs |K_S+B+M-N|$,
where $N=f_*N_1$; i.e., 
$\exists \Lambda \in |K_S+B+M-N|$ with $p \notin \Supp(\Lambda)$. Then
$\Lambda + N \in |K_S+B+M|$ and $p \notin \Supp(\Lambda + N)$, as required.

Note that we haven't used (1.3.3) yet; all we needed so far was $\beta_1 > 0$.

\vspace{8pt}

{\bf (1.9)} Finally, assume that $(S,B,D)$ is not $PLC$ at $p$. Then 
$b_j+d_j > 1$ for some $j$ with $p \in C_j$. In fact, since $2 = \ord_p(B+D) =
\sum (b_i+d_i) \cdot \mult_p(C_i)$, there can be at most one $C_j$ through $p$
with $b_j+d_j > 1$, and then that $C_j$ must be smooth at $p$ and also
$b_i+d_i < 1$ for all $i \neq j$ with $p \in C_i$. Let that $j$ be $0$; thus
$b_0+d_0 > 1$, $C_0$ is smooth at $p$, and $b_i+d_i < 1$ if $i \neq 0$ and
$p \in C_i$. We say that $C_0$ is the {\bf critical curve} at $p$.

  Let $c$ be the {\bf $PLC$ threshold} of $(S,B,D)$ at $p$: 
\begin{equation*}
  c = \max \{ \lambda \in \Q_+ \mid (S,B,\lambda D) \text{ is $PLC$ at $p$} \};
\end{equation*}
explicitly, $b_0+cd_0 = 1$, i.e. $c = \dfrac{1-b_0}{d_0}$. Note that $0<c<1$.

$M-cD \qle (1-ct)M$ is still nef and big on $S$, and we have: 
\begin{equation*}
  K_S+\rup{M-cD} = K_S+B+M-\rdn{B+cD} = K_S+B+M-C_0-N ,
\end{equation*}
with $p \notin \Supp(N)$. (If $p \in C_i$ and $i \neq 0$ then $b_i+d_i < 1$,
and therefore $b_i+cd_i < 1$, too, because $c < 1$; hence $p \notin \Supp(N)$.)

(1.2.1)(a) yields $H^1(S, K_S+B+M-C_0-N) = 0$, and therefore the restriction
map $H^0(S, K_S+B+M-N) \to H^0(C_0, (K_S+B+M-N)|_{C_0})$ is surjective. Hence
it suffices to show that $p \notin \Bs|(K_S+B+M-N)|_{C_0}|$.

\vspace{3pt}

We have
\begin{equation*}
  K_S+B+M-N = K_S+\rup{M-cD}+C_0, 				\tag{1.9.1}
\end{equation*}
  and therefore $(K_S+B+M-N)|_{C_0} = K_{C_0}+\rup{M-cD}|_{C_0}$;
  by (1.2.2)(a), it suffices to show that $\rup{M-cD} \cdot C_0 \geq 2$.
  In any event $\rup{M-cD}\cdot C_0$ is an integer; we will show that
  $\rup{M-cD} \cdot C_0 > 1$.

\vspace{4pt}

$\rup{M-cD} = (M-cD) + \Delta$, where $\Delta = \rup{M-cD}-(M-cD) =
\rup{(M+B)-(B+cD)}-(M-cD) = (M+B)-\rdn{B+cD}-(M-cD) = (B+cD)-\rdn{B+cD}
=\frp{B+cD}$. $\Delta$ is an effective divisor which intersects $C_0$ properly,
because $C_0$ has integer coefficient (namely, 1) in $B+cD$. Moreover, in a 
neighborhood of $p$ we have $\frp{B+cD} = (B+cD) - C_0$, because
$B+cD = C_0 + \sum_{i \neq 0} (b_i+cd_i) C_i$, and $0 \leq b_i+cd_i < 1$
for every $i \neq 0$ such that $p \in C_i$. In particular, we have 
\begin{equation*}
  \ord_p(\Delta) = \ord_p(B+cD)-1 = \mu+c(2-\mu)-1.
\end{equation*}

\vspace{4pt}

\begin{align*}
  \rup{M-cD} \cdot C_0 &= (M-cD) \cdot C_0 + \Delta \cdot C_0 \geq
      (1-ct)M \cdot C_0 + \ord_p(\Delta)	\tag{\bf 1.10}	\\
      &\geq (1-ct)\beta_1 + \mu + c(2-\mu) - 1.
\end{align*}
Therefore the inequality $\rup{M-cD} \cdot C_0 > 1$ follows from 
\begin{equation*}
  (1-ct)\beta_1 > (1-c)(2-\mu).			\tag{1.10.1}
\end{equation*}

If $\beta_1 \geq 2-\mu$ then (1.10.1) is trivial, because
$t < 1 \implies 1-ct > 1-c$.

\vspace{6pt}

{\bf (1.11)} When $\mu < 1$ the inequality we assume for $\beta_1$ (namely, 
(1.3.4)) is weaker than $\beta_1 \geq 2-\mu$. However, in this case the
equation $B+cD = C_0 + \textit{other terms}$ yields a nontrivial lower bound
for $c$: $\mu + c(2-\mu) = \ord_p(B+cD) \geq \ord_p(C_0) =1$, and therefore
$c \geq \dfrac{1-\mu}{2-\mu} > 0$.

The inequality (1.10.1) can also be written as 
\begin{equation*}
  c(2-\mu-t\beta_1) > 2-\mu-\beta_1.		\tag{1.11.1}
\end{equation*}

We may assume that $\beta_1 < 2-\mu$ (or else (1.10.1) is already proved). 
We have $c \geq \dfrac{1-\mu}{2-\mu}$, $t < \dfrac{2-\mu}{\beta_2}$
(see (1.6.1)), and $\dfrac{1-\mu}{\beta_2} \leq 1 - \dfrac{1}{\beta_1}$ 
(by (1.3.4)); therefore
\begin{equation*}
  \begin{split}
    c(2-\mu-t\beta_1) &>\frac{1-\mu}{2-\mu}(2-\mu-\frac{2-\mu}{\beta_2}\beta_1)
  	= (1-\mu-\frac{1-\mu}{\beta_2}\beta_1) 			\\
      &\geq (1-\mu)-(1-\frac{1}{\beta_1})\beta_1 = 2-\mu-\beta_1.
  \end{split}
\end{equation*}

\vspace{3pt}

(1.11.1) is proved. This concludes the proof of Proposition 3.


\section{Separation of points}

{\bf (2.1)} Let $(S,B,M)$ be as in (1.1). Fix two distinct points $p,q\in S$.
In this section we give criteria for $|K_S+B+M|$ to separate $(p,q)$.

Note that in each case $|K_S+B+M|$ is free at $p$ and $q$, by Proposition 3,
and therefore it suffices to find $s \in H^0(S, K_S+B+M)$ such that $s(p)=0,
s(q) \neq 0$, \emph{or} vice-versa.

\begin{Notation}
  $\mu_p = \ord_p(B), \quad \mu_q = \ord_q(B)$.
\end{Notation}

\vspace{6pt}

\begin{Proposition}
  $|K_S+B+M|$ separates $(p,q)$ in each of the following cases:
    \begin{enumerate}
      \item[{\bf 1.}] $\mu_p \geq 2$ and $\mu_q \geq 2$;
      \item[{\bf 2.}] $\mu_q \geq 2$; $0 \leq \mu_p < 2$; $M^2 > (\beta_2)^2,
                      M \cdot C \geq \beta_1$ for every irreducible curve
                      $C \subset S$ passing through $p$, where $\beta_2,
                      \beta_1$ are positive numbers which satisfy (1.3.2)
                      and (1.3.3) for $\mu = \mu_p$;
      \item[{\bf 3.}] $0 \leq \mu_p < 2$ and $0 \leq \mu_q < 2$; 
      		      $M^2 > (\beta_{2,p})^2+(\beta_{2,q})^2$, and 
      		    \begin{enumerate}
      		      \item[(i)] $M \cdot C \geq \beta_{1,p}$ for every
      		          curve $C \subset S$ passing through $p$,
      		      \item[(ii)] $M \cdot C \geq \beta_{1,q}$ for every
      		          curve $C \subset S$ passing through $q$,
      		      \item[(iii)] $M \cdot C \geq \beta_{1,p}+\beta_{1,q}$
      		          if $C$ passes through \emph{both} $p$ \emph{and} $q$,
      		    \end{enumerate}
                      where $\beta_{2,p},\beta_{1,p}\,;\beta_{2,q},\beta_{1,q}$
                      are positive numbers which satisfy the inequalities
    \end{enumerate}
  \begin{align*}
    \beta_{2,p} &\geq 2-\mu_p ,\quad \beta_{2,q} \geq 2-\mu_q;	\tag{2.1.1}  \\
    \beta_{1,p} &\geq \min \left\{ (2-\mu_p) ;
      \frac{\beta_{2,p}}{\beta_{2,p}-(1-\mu_p)} \right\} ,
      \text{ and similarly for $\beta_{1,q}$}.			\tag{2.1.2}
  \end{align*}
\end{Proposition}

\vspace{6pt}

\noindent {\bf Proof of Proposition 4.}

\vspace{4pt}

{\bf (2.2)} Let $f:S_1 \to S$ be the blowing-up of $S$ at $p$ and $q$, with
exceptional curves $E_p, E_q$. As in (1.4), we have: 
\begin{equation*}
  K_{S_1}+\rup{f^*M} = f^*(K_S+B+M) - (\rdn{\mu_p}-1)E_p -
  		(\rdn{\mu_q}-1)E_q.
\end{equation*}

In particular, if $\mu_p \geq 2$ and $\mu_q\geq 2$ (case 1 of the proposition),
we get 
\begin{equation*}
  H^1(S, \Ocal_S(K_S+B+M) \tens {\MM_p}^{t_p} \tens {\MM_q}^{t_q})=0
\end{equation*}
  for positive integers $t_p, t_q$ (compare to (1.5.2)); the conclusion
  follows as in (1.5).

\vspace{8pt}

{\bf (2.3)} Next assume that $\mu_p < 2, \mu_q \geq 2,\, M^2 > (\beta_2)^2$
with $\beta_2 \geq 2-\mu_p,$ etc. (case 2 of the proposition).
Write $\mu = \mu_p$. As in (1.6.1), we can find an effective \Q-divisor $D$
on $S$ such that $\ord_p(D) = 2-\mu$ and $D \qle tM$ for some
$t \in \Q\, , 0 < t < \dfrac{2-\mu}{\beta_2}$.

If $(S,B,D)$ is $PLC$ at $p$, the argument of (1.8) yields a vanishing 
\begin{equation*}
  H^1(S_1, f^*(K_S+B+M) - E_p - N_0) = 0	\tag{2.3.1}
\end{equation*}
where $N_0$ is an effective divisor supported away from $E_p$. Note that in
this case $N_0 \geq E_q$, because $\mu_q \geq 2$. Indeed, (2.3.1) is
obtained by applying (1.2.1)(b) to 
\begin{equation*}
  \begin{split}
    K_{S_1}+\rup{f^*(M-D)} &= f^*(K_S+B+M) - E_p -
            t_q E_q - \sum \rdn{b_i+d_i} D_i   \\
      &= f^*(K_S+B+M) - E_p - t_q E_q - \textstyle\sum ' D_i - N_1,      
  \end{split}					\tag{2.3.2}
\end{equation*}
where $\sum ' D_i$ and $N_1$ are as in (1.8.1) and 
$t_q = \rdn{\mu_q + \ord_q(D)}-1  $ is an integer, $t_q \geq 1$;
then $N_0 = N_1 + t_q E_q \geq E_q$.

\vspace{4pt}

The vanishing (2.3.1) implies the surjectivity of the restriction map 
\begin{multline*}
  H^0(S_1, f^*(K_S+B+M)-N_0)						\\
  	\lra H_0(E_p, (f^*(K_S+B+M)-N_0)|_{E_p}) \cong \C
\end{multline*}
(note that $f^*(K_S+B+M)|_{E_p}$ is trivial, and so is $N_0|_{E_p}$ because 
$N_0 \inters E_p = \emptyset$).

Hence we can find $\Gamma \in |f^*(K_S+B+M)-N_0|$ such that $\Gamma \inters
E_p = \emptyset$. As $\Gamma + N_0 \in |f^*(K_S+B+M)|$, we have
$\Gamma + N_0 = f^*\Lambda$ for some $\Lambda \in |K_S+B+M|$. Moreover, 
$p \notin \Supp(\Lambda)$, because $f^*\Lambda \inters E_p = \emptyset$, but
$q \in \Supp(\Lambda)$, because $f^*\Lambda = \Gamma + N_0 \geq E_q$. Thus
$|K_S+B+M|$ separates $(p,q)$ in this case.

\vspace{6pt}

{\bf (2.4)} Now assume that $(S,B,D)$ is not $PLC$ at $p$. Let $c, C_0$ be
the $PLC$ threshold and the critical curve at $p$, as in \S 1, (1.9)--(1.11).
Let $\phi:S_2 \lra S$ be the blowing-up of $S$ at $q$ (only), with
exceptional curve $F_q$. Let $C_0' \subset S_2$ be the proper transform of
$C_0$ in $S_2$. Let $p' = \phi^{-1}(p)$. We have: 
\[
  K_{S_2}+\rup{\phi^*(M-cD)} = \phi^*(K_S+B+M) - C_0' - N_0, 
\]
where $p' \notin \Supp(N_0)$, as in (1.9), and $N_0 \geq F_q$, as in (2.3).

The argument in (1.9)--(1.11) shows that there exists $\Gamma \in
|\phi^*(K_S+B+M)-N_0|$ with $p' \notin \Supp(\Gamma)$. Now the proof can be
completed as in the last part of (2.3).

\vspace{10pt}

{\bf (2.5)} Finally, consider the case $\mu_p < 2$ and $\mu_q < 2$, with 
$M^2 > (\beta_{2,p})^2 + (\beta_{2,q})^2$, etc. (case 3 of the proposition).

\vspace{4pt}

As in (1.6.1), we can find $G \in |nM|$ with $\ord_p(G) > \beta_{2,p} n$ and
$\ord_q(G) > \beta_{2,q} n$. Let $r = \max \left\{ \dfrac{2-\mu_p}{\ord_p(G)},
\dfrac{2-\mu_q}{\ord_q(G)} \right\},$ and $D=rG$. Then $\ord_p(D) \geq
2-\mu_p$ and $\ord_q(G) \geq 2-\mu_q$, and at least one of the last two
inequalities is an equality. Without loss of generality we may assume that
$\ord_p(D) = 2-\mu_p$ and $m_q \eqdef \ord_q(D) \geq 2-\mu_q$. We have
$D \qle tM$, with 
\begin{equation*}
  0 < t = rn = \frac{2-\mu_p}{\ord_p(G)}\,n < \frac{2-\mu_p}{\beta_{2,p}}
  		\leq 1;				\tag{2.5.1}
\end{equation*}
also, $m_q = \ord_q(D) = r \cdot \ord_q(G) > r \cdot (\beta_{2,q} n) =
t \beta_{2,q}$, and therefore 
\begin{equation*}
  t < \frac{m_q}{\beta_{2,q}}			\tag{2.5.2}
\end{equation*}
(this is the analogue of (2.5.1) at $q$).

\vspace{4pt}

If $(S,B,D)$ is $PLC$ at $p$, then (1.2.1)(b) yields 
\begin{equation*}
  H^1(S_1, f^*(K_S+B+M) - E_p - N_0) =0,	\tag{2.5.3}
\end{equation*}
with $N_0 \inters E_p = \emptyset, N_0 \geq E_q$ (the computation in (2.3.2)
applies unchanged in this situation). In this case we conclude as in (2.3).

\vspace{6pt}

{\bf (2.6)}
Now assume that $(S,B,D)$ is not $PLC$ at $p$. Let $c, C_0$ be the $PLC$ 
threshold and the critical curve at $p$. (1.2.1)(a) yields 
\begin{equation*}
  H^1(S_1, f^*(K_S+B+M) - D_0 - N_0) = 0, \qquad N_0 \inters E_p = \emptyset.
  						\tag{2.6.1}
\end{equation*}

\vspace{4pt}

If $N_0 \inters E_q \neq \emptyset$, we use (2.6.1) to find 
$\Gamma \in |f^*(K_S+B+M)-N_0|$ which does not pass through $\tilde{p} =
D_0 \inters E_p$; the proof is the same as in (1.9)--(1.11). Then the
conclusion follows as in (2.3).

\vspace{4pt}

Assume that $N_0 \inters E_q = \emptyset$. We discuss separately the subcases
$q \in C_0$ and $q \notin C_0$. If $q \in C_0$, we separate $(p,q)$ on $C_0$.
If $q \notin C_0$, we reverse the roles of $p$ and $q$.

\vspace{6pt}

{\bf (2.7)}
First consider the subcase $q \in C_0$. The vanishing (2.6.1) implies 
  \begin{equation*}
      H^1(S, K_S+B+M-C_0-N) = 0,			\tag{2.7.1}
  \end{equation*}
with $N=f_*N_0,\, \Supp(N) \inters \{ p,q \} = \emptyset$. Consequently, the
restriction map 
  \[
      H^0(S, K_S+B+M-N) \lra H^0(C_0, (K_S+B+M-N)|_{C_0})
  \]
is surjective, and it suffices to show that $|(K_S+B+M-N)|_{C_0}|$ separates
$(p,q)$ on $C_0$. As in (1.9.1), we have 
\[
  (K_S+B+M-N)|_{C_0} = K_{C_0} + \rup{M-cD}|_{C_0};
\]
by (1.2.2)(b) it is enough to show that $\rup{M-cD} \cdot C_0 > 2$ (and
consequently $\geq 3$). 

We proceed as in \S 1: $\rup{M-cD}=(M-cD)+\Delta$, with $\Delta=\frp{B+cD}$; 
$\Delta$ and $C_0$ intersect properly, and
$\ord_p(\Delta) = \mu_p+c(2-\mu_p)-1, \ord_q(\Delta) = \mu_q+c m_q-1$
(note that $N_0 \inters E_q = \emptyset \implies$ the only component with
coefficient $\geq 1$ of $B+cD$ through $q$ is $C_0$, and moreover $C_0$ must
be smooth at $q$). Therefore 
\begin{multline*}
    \rup{M-cD} \cdot C_0  = (M-cD) \cdot C_0 + \Delta \cdot C_0	\\
        \geq (1-ct) M \cdot C_0 + \ord_p(\Delta) + \ord_q(\Delta)	\\
        \geq (1-ct)(\beta_{1,p}+\beta_{1,q}) + (\mu_p+c(2-\mu_p)-1) +
        		(\mu_q + c m_q -1)
\end{multline*}
($M \cdot C_0 \geq \beta_{1,p}+\beta_{1,q}$, because this time $C_0$ passes
through both $p$ and $q$.)

Hence $\rup{M-cD} \cdot C_0 > 2$ follows from  
\begin{equation*}
  (1-ct)(\beta_{1,p}+\beta_{1,q})+(\mu_p+c(2-\mu_p)-1)+(\mu_q+c m_q-1) > 2,
  							\tag{2.7.2}
\end{equation*}
which in turn follows from the following two inequalities: 
\begin{align*}
  (1-ct)\beta_{1,p} + (\mu_p+c(2-\mu_p)-1)  &> 1 \quad \text{and}
  				\tag{2.7.3}			\\
  (1-ct)\beta_{1,q} + (\mu_q+c m_q-1)  &> 1.   \tag{2.7.4}
\end{align*}

(2.7.3) is proved like (1.10.1) in \S 1: if $\beta_{1,p} \geq 2-\mu_p$, then
$t < 1 \implies (1-ct) \beta_{1,p} > (1-c)(2-\mu_p) \implies \text{(2.7.3)}$.
If $\beta_{1,p} < 2-\mu_p$ (which can happen only if $\mu_p < 1$), then we
have $c \geq \dfrac{1-\mu_p}{2-\mu_p}$ as in (1.11),
$t < \dfrac{2-\mu_p}{\beta_{2,p}}$ by (2.5.1), and
$\dfrac{1-\mu_p}{\beta_{2,p}} \leq 1 - \dfrac{1}{\beta_{1,p}}$ by (2.1.2), and
therefore (2.7.3) follows as in (1.11).

\vspace{4pt}

(2.7.4) is proved similarly. First, since $m_q = \ord_q(D) \geq 2-\mu_q$, the
inequality is true when $\beta_{1,q} \geq 2-\mu_q$, as in the proof of (2.7.3)
above. When $\beta_{1,q} < 2-\mu_q$ we must have $\mu_q < 1$; then
$B+cD\geq C_0 \implies \mu_q+ cm_q\geq 1 \implies c\geq\dfrac{1-\mu_q}{m_q},
\; t < \dfrac{m_q}{\beta_{2,q}}$ by (2.5.2), and $\dfrac{1-\mu_q}{\beta_{2,q}}
\leq 1-\dfrac{1}{\beta_{2,q}}$ by (2.1.2); consequently 
\begin{equation*}
  \begin{split}
    c(m_q-t\beta_{1,q}) &>
    \frac{1-\mu_q}{m_q} \left( m_q-\frac{m_q}{\beta_{2,q}}\beta_{1,q} \right)
    = 	\\
      &= (1-\mu_q) - \frac{1-\mu_q}{\beta_{2,q}}\beta_{1,q}
      \geq 2-\mu_q-\beta_{1,q},
  \end{split}
\end{equation*}
which yields (2.7.4).

Thus (2.7.2) is proved; this concludes the proof when $q \in C_0$.

\vspace{6pt}

{\bf (2.8)} To complete the proof of the proposition in case 3, consider the
remaining subcase, $q \notin C_0$. In this subcase separation of $(p,q)$ is
obtained by reversing the roles of $p$ and $q$. Namely, let $D'=\alpha D$, for
the positive rational number $\alpha$ such that $\ord_q(D')=2-\mu_q$; that is,
$\alpha = \dfrac{2-\mu_q}{\ord_q(D)} = \dfrac{2-\mu_q}{m_q}$. Note that
$D' \qle t'M$, where $t' = \alpha t <
\dfrac{2-\mu_q}{m_q} \cdot \dfrac{m_q}{\beta_{2,q}}$ (by (2.5.2)), i.e. 
\begin{equation*}
  0 < t' < \frac{2-\mu_q}{\beta_{2,q}} \leq 1.	\tag{2.8.1}
\end{equation*}

Let $c'$ be the $PLC$ threshold for $(S,B,D')$ at $q$; note that $c'\alpha
> c$ ($c'\alpha$ is the $PLC$ threshold of $(S,B,D)$ at $q$, and therefore
$c < c'\alpha$ follows from $N_0 \inters E_q = \emptyset$ in (2.6.1)).
This, in turn, implies $B+c'D' = B+c'\alpha D \geq C_0$.

If $(S,B,D')$ is $PLC$ at $q$ (i.e. if $c'=1$), then (1.2.1)(b) yields 
\begin{equation*}
  H^1(S_1, f^*(K_S+B+M)-E_q-N_0') = 0, \qquad N_0' \inters E_q = \emptyset
  						\tag{2.5.$3'$}
\end{equation*}
(Compare to (2.5.3)).

If $(S,B,D')$ is not $PLC$ at $q$ (i.e. if $c' < 1$), and $C_0'$ is the
critical curve at $q$, then (1.2.1.)(a) yields 
\begin{equation*}
  H^1(S_1, f^*(K_S+B+M)-D_0'-N_0') = 0,	\qquad N_0' \inters E_q = \emptyset
  						\tag{2.6.$1'$}
\end{equation*}
(Compare to (2.6.1), noting that now $p$ and $q$ are interchanged.)

In both cases, the arguments in (1.8) and, respectively, (1.9)--(1.11) show
that there exists $\Lambda \in |K_S+B+M-N'|$ with $q \notin \Supp(\Lambda)$,
where $N' = f_*N_0'$ is an effective divisor with $q \notin \Supp(N')$. Now,
however, $N' \geq C_0$ (because $B+c'D' \geq C_0$, as noted earlier, and
$q \notin C_0 \implies C_0$ is not discarded even when the vanishing theorem
is used in the form (1.2.1)(b)); thus $\Gamma + N' \in |K_S+B+M|$ passes
through $p$ but not through $q$.

\vspace{6pt}

This completes the proof of Proposition 4.


\section{Separation of tangent directions}

{\bf (3.1)} Let $(S,B,M)$ be as in \S 1. Fix a point $p \in S$. In this section
we give criteria for $|K_S+B+M|$ to separate directions at $p$.

The statements (and proofs) are somewhat similar to those in \S 2. The main 
difference is in the part of the proof corresponding to the discussion in
(2.8).  So far in our proofs we worked with $M-cD$, where $c$ was always the 
$PLC$ threshold at some point or another; this made the arguments relatively
transparent. In (2.8), when we passed from $c$ = $PLC$ threshold at $p$ to
$c'\alpha$ = $PLC$ threshold at $q$, the relevant fact was that $q \notin C_0$,
where $C_0$ was the critical curve at $p$, and therefore $C_0$ did not affect
the local computations around $q$. In separating tangent directions, the 
analogue is a curve $C_0$ through $p$, such that $\vv \notin T_p(C_0)$ for some
fixed $\vv \in T_p(S), \vv \neq \vec{0}$. Then we will have to increase $c$
to some larger value $c'$, but clearly in that case $(S, B, c'D)$ will no 
longer be $PLC$ at $p$. While this complicates the computations, the geometric
idea is still the same: find a divisor $\Gamma \in |K_S+B+M-C_0-N|$, 
$p \notin \Supp(N)$, such that $\Gamma$ does not pass through $p$; then
$\Gamma + C_0 +N$ has only one component through $p$, namely, $C_0$, and
$\vv \notin T_p(C_0)$ -- therefore $\Gamma + C_0 +N$ passes through $p$ and
is not tangent to $\vv$, as required.

Another technical problem, which did not arise before, is that in some cases
the ``minimizing'' curve $C_0$ may be singular at $p$. (This possibility is
directly related to the need, in some cases, to increase $c$ beyond the 
$PLC$ threshold at $p$.) In those cases we separate the tangent directions
on $C_0$, using (1.2.2)(b) (note that $C_0$ singular at $p$ \ $\implies 
T_p(S)=T_p(C_0)$); the vanishing (1.2.1) is then used to lift from $C_0$
to $S$.

\vspace{6pt}

{\bf (3.2)} Let $S$ be a smooth surface, as before; let $p$ denote a point on
$S$, and fix $\vv \in T_p(S), \vv \neq \vec{0}$. Let $Z$ denote the 
zero-dimensional subscheme of length $2$ of $S$, corresponding to $(p,\vv)$;
in local coordinates $(x,y)$ at $p$ such that $\vv$ is tangent to $(y=0)$,
$Z$ is defined by the ideal $\I_Z = (x^2,y) \cdot \Ocal_S$.

Let $f:S_1 \to S$ be the blowing-up of $S$ at $p$, with exceptional curve
$E_p$, and let $V \in E_p$ correspond to (the direction of) $\vv$. Let
$g:S_2 \to S_1$ be the blowing-up of $S_1$ at $V$, with exceptional curve
$F_{\vv}$, and let $F_p = g^{-1}E_p$. Let $h = g \circ f$. Write 
\begin{equation*}
  h^*B = h^{-1}B + \mu_p F_p + \mu_{\vv} F_{\vv};	\tag{3.2.1}
\end{equation*}
$\mu_p = \ord_p(B)$, while (3.2.1) is the definition of $\mu_{\vv}$.

More generally, if $G$ is any effective \Q-divisor on $S$, denote the order of 
$h^*G$ along $F_{\vv}$ by $o_{\vv}(G)$; $o_{\vv}(G) = \ord_p(G) +
\ord_V(f^{-1}G)$. For convenience, let $o_V(G) \eqdef \ord_V(f^{-1}G)$, and
let $\mu_V = o_V(B)$.

Note that, in general, $o_{\vv} = \ord_p + o_V$ and $o_V \leq \ord_p$; in
particular: 
\begin{equation*}
  \mu_p \leq \mu_{\vv} \leq 2\mu_p.			\tag{3.2.2}
\end{equation*}

\vspace{6pt}

{\bf (3.3)} Consider again $(S,B,M)$ as in \S 1, and fix $p, \vv$ as in (3.2).
Since the proofs will now be more complex, we will state the criteria for
separating $\vv$ at $p$ one by one, in increasing order of difficulty.

The first (and easiest) case is:

\begin{Proposition}[Case 1]
  If $\mu_p \geq 3$ or $\mu_{\vv} \geq 4$, then $|K_S+B+M|$ separates $\vv$
  at $p$. ($M$ must still be nef and big.)
\end{Proposition}

\begin{proof}
    Recall that the conclusion means that the restriction map
  $$
    H^0(S, K_S+B+M) \to H^0(Z, K_S+B+M|_Z) \cong \Ocal_Z
  $$
  is surjective.

    If $\mu_p \geq 3$, we use the vanishing theorem in the form (1.2.1)(a) for
  \begin{align*}
    K_{S_1} + \rup{f^*M} &= f^*(K_S+B+M) + E_p - \rdn{f^*B}	\\
                         &= f^*(K_S+B+M) -t E_p,
  \end{align*}
  where $t = \rdn{\mu_p} - 1 \geq 2$, as in (1.4)--(1.5); then
  $H^0(S, K_S+B+M) \to \Ocal_S/{\MM_p}^t$ is surjective, and since $t \geq 2$,
  we have ${\MM_p}^t \subset \I_Z$, i.e. $\Ocal_S/{\MM_p}^t \to \Ocal_Z$ is
  also surjective. (See also Remark 1.5.3.)

  If $\mu_{\vv} \geq 4$, the agrument is similar, starting on $S_2$:
  \begin{align*}
    K_{S_2}+\rup{h^*M} &= h^*(K_S+B+M)+F_p+2F_{\vv}-\rdn{h^*B}	\\
    		       &= h^*(K_S+B+M)-t_pF_p - t_{\vv}F_{\vv} ,
  \end{align*}
  where $t_{\vv} = \rdn{\mu_{\vv}}-2 \geq 2$, and $t_p = \rdn{\mu_p}-1 \geq 1$
  (indeed, by (3.2.2), $\mu_p \geq \frac{1}{2}\mu_{\vv} \geq 2$.) As in the
  previous case, we get a vanishing $H^1(S, \Ocal_S(K_S+B+M) \tens \I) = 0$ for
  $\I = h_*\Ocal_{S_2}(-t_pF_p-t_{\vv}F_{\vv})$; $\Supp(\Ocal_S/\I) =
  \{ p \}$ and $\I \subset \I_Z$, so the conclusion follows as before.
\end{proof}

\vspace{6pt}

{\bf (3.4)} Now assume that $\mu_p < 3$ and $\mu_{\vv} < 4$. 

First consider the case $2 \leq \mu_p < 3$. Then $2 \leq \mu_{\vv} <4$, and
therefore $0 < (4-\mu_{\vv}) \leq 2$.

\addtocounter{Theorem}{-1}
\begin{Proposition}[Case 2]
  Let $2 \leq \mu_p < 3$ and $2 \leq \mu_{\vv} < 4$. Assume that $M^2 > 
  (4-\mu_{\vv})^2$, $M \cdot C \geq \frac{1}{2}(4-\mu_{\vv})$ for every curve
  $C \subset S$ through $p$, and $M \cdot C \geq (4-\mu_{\vv})$ for every
  curve $C$ containing $Z$ -- i.e., such that $p \in C$ and $\vv \in T_p(C)$.
  Then $|K_S+B+M|$ separates $\vv$ at $p$.
\end{Proposition}

\vspace{4pt}

\emph{Proof.}

{\bf (3.5) Claim:} We can find an effective \Q-divisor $D$ on $S$ such that
  $o_{\vv}(D) = 4-\mu_{\vv}$ and $D \qle tM$ for some $t \in \Q\,, 0 < t < 1$.
  (See (3.2) for the definition of $o_{\vv}(D)$.)

\emph{Proof of (3.5):} Choose $a > (4-\mu_{\vv})$ such that $M^2 > a^2$. Then
  $(h^*M-a F_{\vv})^2 = M^2-a^2 >0$ and $(h^*M-a F_{\vv}) \cdot h^*M=M^2>0$;
  therefore $h^*M-a F_{\vv} \in N(S_2)^+$, the positive cone of $S_2$, and in
  particular it is big. (See, for example, \cite[(1.1)]{km}.) Therefore
  $\exists T$, effective \Q-divisor on $S_2$, such that $T \qle h^*M -
  a F_{\vv}$. Put $D_1 = h_*(T+a F_{\vv})$; then $D_1 \qle h_*(h^*M)=M$. Also,
  $h^*D_1= T + a F_{\vv}$ (their difference has support contained in $F_p
  \union F_{\vv}$; on the other hand, $T + a F_{\vv} \qle h^*M \implies
  \big( h^*D_1 - (T+a F_{\vv}) \big) \cdot F_p = 
  \big( h^*D_1 - (T+a F_{\vv}) \big) \cdot F_{\vv} = 0$, and 
  $h^*D_1 = T+a F_{\vv}$ follows from the negative definiteness of the
  intersection form on $h^{-1}(p) = F_p \union F_{\vv}$).
  We have $D_1 \qle M$ and $o_{\vv}(D_1) \geq a > 4 - \mu_{\vv}$. Take
  $D = tD_1$, $t = \dfrac{4-\mu_{\vv}}{o_{\vv}(D_1)} < 1$.    \qed

\vspace{4pt}

\begin{Remark}
  The statement of (3.5) is similar to that of (1.6.1), and indeed, we could
  have proved it as in \S 1. However, the proof we gave here is easier to
  generalize, especially on \emph{singular} surfaces.
\end{Remark}

\vspace{4pt}

{\bf (3.6)} We return to the proof of Proposition 5, Case 2. Choose $D$ as in
(3.5). Write $B=\sum b_iC_i, D=\sum d_iC_i; D_i=f^{-1}C_i, T_i=g^{-1}D_i =
h^{-1}C_i; h^*B=h^{-1}B + \mu_p F_p + \mu_{\vv} F_{\vv}, h^*D = h^{-1}D +
m_p F_p + (4-\mu_{\vv})F_{\vv}$, where $m_p=\ord_p(D)$. 
We have $K_{S_2} = h^*K_S +F_p +2F_{\vv}$.

If $b_i+d_i \leq 1$ for every $C_i$ through $p$, then 
\begin{align*}
  K_{S_2}+\rup{h^*(M-cD)} &=h^*(K_S+B+M)+F_p+2F_{\vv}-\rdn{h^*(B+D)}
  				\\
     &= h^*(K_S+B+M)-t_pF_p- 2 F_{\vv}-\textstyle\sum'T_i-N_2,  \tag{3.6.1}
\end{align*}
where $\sum ' T_i$ extends over all $i$ with $b_i+d_i=1$ and $p \in C_i$ (if
any), $N_2$ is an effective divisor on $S_2$ such that $\Supp(N_2) \inters
h^{-1}(p) = \emptyset$, and $t_p = \rdn{\mu_p + m_p}-1 \geq 1$ (because 
$\mu_p \geq 2$ by hypothesis). Then we conclude as in (3.3) (Case 1 of the
Proposition), using the vanishing (1.2.1)(b) to dispose of $\sum ' T_i$ (if it
is not zero).

\vspace{4pt}

{\bf (3.7)} Now assume that $b_i+d_i > 1$ for at least one $C_i$ through $p$.
Let 
\begin{equation*}
  c \eqdef \min \left\{ \frac{3-\mu_p}{m_p} ; \frac{1-b_i}{d_i} :
  	b_i+d_i>1 \text{ and } p \in C_i \right\}.		\tag{3.7.1}
\end{equation*}

\vspace{4pt}

If $c = \dfrac{3-\mu_p}{m_p}$, we finish again as in Case 1, using (1.2.1)(b) 
for
\[
  K_{S_1} + \rup{f^*(M-cD)} = f^*(K_S+B+M) -2E_p - \textstyle\sum ' D_i - N_1
\]
on $S_1$, where $\sum ' D_i$ extends over all $i$ such that $b_i+cd_i=1$ and
$p \in C_i$ (if any), and $\Supp(N_1) \inters E_p = \emptyset$.

\vspace{4pt}

{\bf (3.8)} If $c = \dfrac{1-b_0}{d_0} < \dfrac{3-\mu_p}{m_p}$ for some $C_0$
through $p$, then
\[
  \sum (b_i+cd_i) \cdot \mult_p(C_i) = \mu_p + cm_p <3 ;
\]
therefore $\mult_p(C_0) \leq 2$, and moreover, if $\mult_p(C_0) = 2$, then
$b_i+cd_i < 1$ for all $C_i$ through $p$ with $i \neq 0$. 
Also, $\mu_{\vv}+c(4-\mu_{\vv}) < 4$ (since $c<1$), and therefore
$B+cD \geq C_0 \implies o_{\vv}(C_0) \leq 3$.

\vspace{4pt}

{\bf (3.9)} If $C_0$ is singular at $p$ and $\vv \notin TC_p(C_0)$ (the 
\emph{tangent cone} to $C_0$ at $p$), then $o_{\vv}(C_0) = 2$. 
We have $Z \subset C_0$, and 
\begin{align*}
  K_S + \rup{M-cD} &= (K_S+B+M) - \rdn{B+cD}				\\
  		   &= (K_S+B+M) - C_0 - N ,
\end{align*}
with $p \notin \Supp(N)$. Using (1.2.1)(a), as in \S 1, it suffices to show
that $\big( (K_S+B+M)-C_0-N \big) |_{C_0}$ separates $\vv$ at $p$ on $C_0$;
that, in turn, will follow from (1.2.2)(b), \emph{if} we can show that 
$\rup{M-cD} \cdot C_0 > 2$.

As before, write $\rup{M-cD} = (M-cD)+\Delta$; $\Delta = \frp{B+cD}$ and $C_0$
intersect properly, and $\Delta = B+cD-C_0$ in an open neighborhood of $p$. 

We have: $\ord_p(\Delta) = \mu_p + cm_p - 2$, and therefore $\Delta \cdot C_0
\geq 2(\mu_p +cm_p -2)$. However, we get a better estimate if we consider
orders along $F_{\vv}$, as follows:
  $o_{\vv}(\Delta) = \mu_{\vv} + c(4-\mu_{\vv}) - 2, \quad \text{because }
  	o_{\vv}(C_0) = 2 $;
$\ord_p(\Delta) \geq \frac{1}{2} o_{\vv}(\Delta)$, and therefore
\[
  \Delta \cdot C_0 \geq \frac{1}{2} o_{\vv}(\Delta) \cdot 2 \geq
  		\mu_{\vv} + c(4-\mu_{\vv}) - 2 .
\]
Finally, 
\begin{align*}
    \rup{M-cD} \cdot C_0 &= (M-cD) \cdot C_0 + \Delta \cdot C_0
              = (1-ct) M \cdot C_0 + \Delta \cdot C_0 \\
    &\geq (1-ct)(4-\mu_{\vv}) + \mu_{\vv} + c(4-\mu_{\vv}) - 2 \tag{3.9.1}  \\
    &> 2 \qquad \text{(because $t<1$),}
\end{align*}
as required

\vspace{4pt}

{\bf (3.10)} If $C_0$ is singular at $p$ and $\vv \in TC_p(C_0)$, then
$o_{\vv}(C_0) = 3$ ($ \geq 3$ is clear, and $\leq 3$ was shown in (3.8)).

Working as in (3.9), we can show that 
\begin{equation*}
  \rup{M-cD} \cdot C_0 \geq (1-ct)(4-\mu_{\vv})+\mu_{\vv}+c(4-\mu_{\vv})-3
  	> 1							\tag{3.10.1}
\end{equation*}
(now $o_{\vv}(\Delta) = o_{\vv}(B+cD-C_0) = \mu_{\vv}+c(4-\mu_{\vv})-3$); thus
in this case we cannot use (1.2.2)(b) as in (3.9). We will modify the argument
as follows:

Start with $f^*(M-cD)$ on $S_1$; the vanishing theorem yields 
\begin{equation*}
  H^1(S_1, f^*(K_S+B+M)-E_p-D_0-N_1)=0, \qquad N_1 \inters E_p = \emptyset
  								\tag{3.10.2}
\end{equation*}
(the coefficient of $E_p$ is $-1$ because $2 \leq \mu_p+cm_p < 3$; the first
inequality follows from $\mu_p \geq 2$, and the second was shown in (3.8)).

$\vv \in TC_p(C_0) \implies V \in D_0$ (recall that $V \in E_p$ corresponds to
$\vv \in T_p(S)$). (3.10.2) implies the surjectivity of the restriction map
\begin{multline*}
  H^0(S_1, f^*(K_S+B+M)-E_p-N_1) 				\\
  		\to H^0(D_0, f^*(K_S+B+M)-E_p-N_1 |_{D_0}).
  								\tag{3.10.3}
\end{multline*}
We will show that $\exists \tilde{\Gamma} \in \left| f^*(K_S+B+M)-E_p-N_1
|_{D_0} \right|$ such that $V \notin \Supp(\tilde{\Gamma})$. Then we can lift
$\tilde{\Gamma}$ to $\Gamma \in |f^*(K_S+B+M)-E_p-N_1|$, since (3.10.3) is
surjective. $\Gamma+E_p+N_1 \in |f^*(K_S+B+M)|$ has the form $f^*\Lambda$ for
some $\Lambda \in |K_S+B+M|$. Finally, $p \in \Supp(\Lambda)$, because
$f^*\Lambda \geq E_p$, but $\vv \notin T_p(\Lambda)$, because
$V \notin \Supp(f^*\Lambda - E_p)$; this shows that $|K_S+B+M|$ separates $\vv$
at $p$ on $S$.

To prove the existence of $\tilde{\Gamma}$, note that
$\big( f^*(K_S+B+M)-E_p-N_1 \big) |_{D_0} = K_{D_0} + \rup{f^*(M-cD)}|_{D_0}$;
we will show that $\rup{f^*(M-cD)} \cdot D_0 > 1$ -- then (1.2.2)(a) implies
the existence of $\tilde{\Gamma}$.

As in (1.9), we can write $\rup{f^*(M-cD)} = f^*(M-cD) + \Delta_1$, where
$\Delta_1 = \frp{f^*(B+cD)}$ and $D_0$ intersect properly, and
$\Delta_1 = f^*(B+cD) - 2E_p - D_0 = f^*(B+cD-C_0)$ in a neighborhood of $E_p$
(the coefficient of $E_p$ in $f^*(B+cD)$ is $\mu_p+cm_p$, and 
$2 \leq \mu_p+cm_p <3$). We have:
\begin{equation*}
  \begin{split}
    \rup{f^*(M-cD)} \cdot D_0 &= f^*(M-cD) \cdot D_0 + \Delta_1 \cdot D_0 \\
    	&\geq (M-cD) \cdot C_0 + \ord_p(B+cD-C_0) \cdot \mult_p(C_0)	  \\
    	&\geq (1-ct)(4-\mu_{\vv})+\mu_{\vv}+c(4-\mu_{\vv})-3		  \\
    	&> 1
  \end{split}							\tag{3.10.4}
\end{equation*}
as in (3.10.1)

\vspace{4pt}

{\bf (3.11)} Now consider the case: $C_0$ smooth at $p$ and tangent to $\vv$,
and $b_i+cd_i < 1$ for all $i \neq 0$ with $p \in C_i$.

  Write $\rup{M-cD} = (M-cD) + \Delta$, where $\Delta$ and $C_0$ intersect
properly; then $\Delta \cdot C_0 = h^*\Delta \cdot T_0$ (projection formula:
recall that $T_0=h^{-1}C_0$) $\geq o_{\vv}(\Delta)$, because 
$F_{\vv} \cdot T_0 = 1$; since $\Delta = \frp{B+cD} = B+cD-C_0$ in a
neighborhood of $p$, we have $o_{\vv}(\Delta) = \mu_{\vv}+c(4-\mu_{\vv})-2$,
and therefore
\begin{equation*}
  \rup{M-cD} \cdot C_0 \geq (1-ct)(4-\mu_{\vv})+\mu_{\vv}+c(4+\mu_{\vv})-2 >2,
\end{equation*}
exactly as in (3.9.1). Thus $K_{C_0}+\rup{M-cD}|_{C_0}$ separates $\vv$ on
$C_0$; we conclude as in (3.9).

\vspace{4pt}

{\bf (3.12)} If $C_0$ is smooth at $p$ and tangent to $\vv$, and moreover
$b_i+cd_i = 1$ for at least one $i \neq 0$ with $p \in C_i$, then: such an $i$
is unique, say $i=1$, and $C_1$ must be smooth at $p$ and not tangent to $\vv$;
indeed, $B+cD \geq C_0+C_1$, while $\ord_p(B+cD) < 3$ and 
$o_{\vv}(B+cD) < 4$.

  In this case reverse the roles of $C_0$ and $C_1$: thus $C_0$ will be smooth
at $p$ and not tangent to $\vv$. This situation is covered below, in (3.13).

\vspace{4pt}

{\bf (3.13)} Finally, assume that $C_0$ is smooth at $p$ and not tangent to
$\vv$. In this case we work with $M-c'D$ for some $c' \geq c$, namely: 
\[
  c' \eqdef \min \left\{ 1; \frac{3-\mu_p}{m_p}; \frac{2-b_0}{d_0};
  	\frac{1-b_i}{d_i} \text{ with } i \neq 0, p \in C_i \text{ and }
  	b_i+d_i \geq 1 \right\}.
\]

In all cases, $M-c'D \qle (1-c't)M$ is still nef and big; using the vanishing
$H^1(S, K_S+\rup{M-c'D})=0$, or the corresponding vanishing on $S_1$ or $S_2$,
we will show that $\exists \Lambda \in |K_S+B+M-C_0-N|, p \notin \Supp(N)$,
such that $p \notin \Supp(\Lambda)$. Then $\Lambda+C_0+N \in |K_S+B+M|$ has
the unique component $C_0$ through $p$ not tangent to $\vv$, as required.

\vspace{4pt}

It remains to prove the existence of $\Lambda$.

\vspace{4pt}

{\bf (3.13.1)} If $c' = \dfrac{3-\mu_p}{m_p}$, then (1.2.1)(b) yields
\[
  H^1(S_1, f^*(K_S+B+M) - 2E_p - D_0 - N_1) = 0, \qquad
  N_1 \inters E_p = \emptyset;
\]
thus $H^1(S_1, f^*(K_S+B+M-C_0-N) - E_p) = 0$, where $N = f_*N_1$, and the
existence of $\Lambda$ follows.

\vspace{4pt}

When $c' = 1$ the proof is similar, starting on $S_2$, as in the proof of
Case~1 of the proposition.

\vspace{4pt}

{\bf (3.13.2)} If $c' = \dfrac{1-b_1}{d_1} < \dfrac{3-\mu_p}{m_p}$ for another
curve $C_1$ through $p$ with $b_1+d_1 > 1$, then $C_1$ must be smooth at $p$
(because $B+c'D \geq C_0+C_1$, and $\ord_p(B+c'D) = \mu_p+c'm_p < 3$). We may
have $c' = c$ (e.g., in the case discussed in (3.12)), or $c' > c$. In any
event, $b_i+c'd_i < 1$ for all curves $C_i$ through $p$, $i \neq 0,1$.
(1.2.1)(a) yields: 
\begin{equation*}
  H^1(S, K_S+B+M-C_0-C_1-N) = 0, \qquad p \notin \Supp(N).	\tag{3.13.3}
\end{equation*}
We claim that $p \notin \Bs \left| K_S+B+M-C_0-N |_{C_1} \right|$, which in
turn follows from (1.2.2)(a) once we show that $\rup{M-c'D} \cdot C_1 > 1$.
Then we use (3.13.3) to lift from $C_1$ to $S$, proving the existence of
$\Lambda$ as stated.

\vspace{3pt}

$\rup{M-c'D} = (M-c'D) + \Delta'$, with $\Delta' =\frp{B+c'D} = B+c'D-C_0-C_1$
in a neighborhood of $p$, and $\Delta', C_1$ intersect properly.

\vspace{3pt}

If $\vv \notin T_p(C_1)$, then $M \cdot C_1 \geq \frac{1}{2}(4-\mu_{\vv})$ by
hypothesis, and $\ord_p(\Delta') \geq \frac{1}{2}o_{\vv}(\Delta') \geq
\frac{1}{2} \big( \mu_{\vv}+c'(4-\mu_{\vv})-2 \big)$; therefore
\[
  \rup{M-c'D} \cdot C_1 \geq \dfrac{1}{2}(1-c't)(4-\mu_{\vv}) +
    \dfrac{1}{2} \big( \mu_{\vv}+c'(4-\mu_{\vv}) - 2 \big) > 1,
\]
as required (compare to (3.9.1)). The proof is the same when
$c' = \dfrac{2-b_0}{d_0}$, i.e. $C_1 = C_0$; in that case
$p \notin \Bs \left| K_S+B+M-C_0-N |_{C_0} \right|$.

\vspace{3pt}

If $\vv \in T_p(C_1)$, then $M \cdot C_1 \geq 4-\mu_{\vv}$ and 
$\Delta' \cdot C_1 \geq o_{\vv}(\Delta') = \mu_{\vv}+c'(4-\mu_{\vv})-3$;
all told, we have
\[
  \rup{M-c'D} \cdot C_1 \geq (1-c't)(4-\mu_{\vv})+\mu_{\vv}+c'(4-\mu_{\vv})-3
  		>1,
\]
as claimed.

\vspace{6pt}

This concludes the proof of Proposition 5, Case 2.

\vspace{10pt}

{\bf (3.14)} Finally, consider the case $0 \leq \mu_p < 2$ (and therefore
$0 \leq \mu_V < 2$ and $0 \leq \mu_{\vv} = \mu_p + \mu_V < 4$).

\vspace{4pt}

\addtocounter{Theorem}{-1}
\begin{Proposition}[Case 3]
  Assume that $0 \leq \mu_p < 2$. Assume, moreover, that 
  $M^2 > (\beta_{2,p})^2 + (\beta_{2,V})^2$ and 
  \begin{enumerate}
    \item[(i)] $M \cdot C \geq \beta_1$ for every curve $C \subset S$ passing
                             through $p$,
    \item[(ii)] $M \cdot C \geq 2\beta_1$ for every curve $C$ containing $Z$
        (i.e., passing through $p$ and with $\vv \in T_p(C)$),
  \end{enumerate}
  where $\beta_{2,p}, \beta_{2,V}, \beta_1$ are positive numbers which
  satisfy:
  \begin{align*}
    \beta_{2,p} &\geq 2-\mu_p, \qquad \beta_{2,V} \geq 2-\mu_V;	\tag{3.14.1} \\
    \beta_1 &\geq \min \left\{ \frac{1}{2}(4-\mu_{\vv});
      \frac{\beta_{2,p}+\beta_{2,V}}{\beta_{2,p}+\beta_{2,V}-(2-\mu_{\vv})}
      \right\}							\tag{3.14.2}
  \end{align*}
\end{Proposition}
 
\begin{proof} 
  The proof is very similar, in many respects, to that of Case 2. We indicate
the main steps of the proof, and we provide explicit computations in a few
cases, to show what kind of alterations are needed.

\vspace{6pt}

{\bf (3.15) Claim:} We can find $D$, an effective \Q-divisor on $S$, such that
$o_{\vv}(D) = 4-\mu_{\vv}$ and $D \qle tM$ for some $t \in \Q\,, t > 0$,
satisfying 
\begin{equation*}
  t < \frac{4-\mu_{\vv}}{\beta_{2,p}+\beta_{2,V}}		\tag{3.15.1}
\end{equation*}
-- and therefore, in particular, $t < 1$.

\vspace{3pt}

\emph{Proof of (3.15):} Choose $a > \beta_{2,p} , b > \beta_{2,V}$, such that
$M^2 > a^2+b^2$. We have $\big( aF_p + (a+b)F_{\vv} \big)^2 = -(a^2+b^2)$, and
therefore $h^*M-\big( aF_p + (a+b)F_{\vv} \big)$ is big, as in the proof of
(3.5). Thus we can find $D_1 \qle M$ on $S$, $D_1 \geq 0$, such that
$h^*D_1 \geq aF_p + (a+b)F_{\vv}$. Then take $D = tD_1$, with 
\[
  t = \frac{4-\mu_{\vv}}{o_{\vv}(D_1)} \leq \frac{4-\mu_{\vv}}{a+b}
    < \frac{4-\mu_{\vv}}{\beta_{2,p}+\beta_{2,V}}. 
\]			\qed

\vspace{6pt}

{\bf (3.16)} If $D = \sum d_i C_i$, as before, and $b_i+d_i \leq 1$ for 
every $C_i$ through $p$, we conclude as in (3.6).

If $b_i+d_i > 1$ for at least one $C_i$ through $p$, then define 
\begin{equation*}
  c = \min \left\{ \frac{3-\mu_p}{m_p}; \frac{1-b_i}{d_i} : 
  	b_i+d_i > 1 \text{ and } p \in C_i \right\}.		\tag{3.16.1}
\end{equation*}

\vspace{4pt}

If $c=\dfrac{3-\mu_p}{m_p}$, we finish as in (3.7).

\vspace{3pt}

If $c = \dfrac{1-b_0}{d_0} < \dfrac{3-\mu_p}{m_p}$ for some $C_0$ through $p$,
then $\mult_p(C_0) \leq 2$ and $o_{\vv}(C_0) \leq 3$; if $C_0$ is singular at
$p$, then it is the only $C_i$ through $p$ with $b_i+cd_i \geq 1$, and we
proceed as in (3.9) or (3.10), according to whether $\vv \in TC_p(C_0)$ or not.
Only the proof of $\rup{M-cD} \cdot C_0 > 2$ (if $\vv \notin TC_p(C_0)$) or
$ > 1$ (if $\vv \in TC_p(C_0)$) needs adjustment.

\vspace{3pt}

Assume first that $\vv \notin TC_p(C_0)$ (with $C_0$ singular at $p$).
Then $\rup{M-cD} = (M-cD) + \Delta, \Delta = \frp{B+cD} = B+cD-C_0$ in a
neighborhood of $p$, and $o_{\vv}(\Delta) = \mu_{\vv}+c(4-\mu_{\vv})-2$;
$\ord_p(\Delta) \geq \frac{1}{2} o_{\vv}(\Delta)$ and $\mult_p(C_0)=2$, so that
\begin{equation*}
  \rup{M-cD} \geq (1-ct)(2\beta_1) +\mu_{\vv}+c(4-\mu_{\vv})-2.  \tag{3.16.2}
\end{equation*}

\vspace{3pt}

If $\beta_1 \geq \frac{1}{2}(4-\mu_{\vv})$, then $\rup{M-cD} \cdot C_0 > 2$
follows from (3.16.2) and $t < 1$. In particular, this is true if 
$\mu_{\vv} \geq 2$. If $\mu_{\vv} < 2$, the hypothesis is weaker than 
$\beta_1 \geq \frac{1}{2}(4-\mu_{\vv})$, namely: 
\begin{equation*}
  \beta_1 \geq
  \frac{\beta_{2,p}+\beta_{2,V}}{\beta_{2,p}+\beta_{2,V} - (2-\mu_{\vv})}.
  								\tag{3.16.3}
\end{equation*}
Assume also that $\beta_1 < \frac{1}{2}(4-\mu_{\vv})$ (otherwise we are done).
Then $\rup{M-cD} \cdot C_0 > 2$ follows from (3.16.2), (3.16.3), (3.15.1), and
$c \geq \dfrac{2-\mu_{\vv}}{4-\mu_{\vv}}$, exactly as in (1.11).

\vspace{4pt}

Now consider the case $\vv \in TC_p(C_0)$ (with $C_0$ still singular at $p$).
Using the strategy of (3.10), all we need to prove is
$\rup{M-cD} \cdot C_0 > 1$, which follows from 
\begin{equation*}
  \rup{M-cD} \cdot C_0 \geq (1-ct)(2\beta_1) +\mu_{\vv}+c(4-\mu_{\vv})-3
  								\tag{3.16.4}
\end{equation*}
(same computation as in (3.10) -- see (3.10.4)). Using (3.16.4), the inequality 
$\rup{M-cD} \cdot C_0 > 1$ is proved exactly as in the previous paragraph.

\vspace{6pt}

{\bf (3.17)} If $C_0$ is smooth at $p$, $\vv \in T_p(C_0)$, and $b_i+cd_i<1$
for every $C_i$ through $p$ with $i \neq 0$, then the proof goes as in (3.11);
the inequality we need in this case, $\rup{M-cD} \cdot C_0 > 2$, is proved as
above.

As in the proof of Case 2 of the Proposition, if $C_0$ is smooth at $p$ and
tangent to $\vv$, and $b_1+cd_1 = 1$ for one more curve $C_1$ through $p$,
then $C_1$ is unique with these properties, and is smooth at $p$ and 
$\vv \notin T_p(C_1)$. Switching $C_0$ and $C_1$, we are in the situation
discussed below. (Compare to (3.12).)

\vspace{6pt}

{\bf (3.18)} Finally, assume that $C_0$ is smooth at $p$ and
$\vv \notin T_p(C_0)$. Define
\[
  c' = \min \left\{ 1; \frac{3-\mu_p}{m_p} ; \frac{2-b_0}{d_0} ;
  \frac{1-b_i}{d_i}: i \neq 0, b_i+d_i > 1 \text{ and } p \in C_i \right\}.
\]
Consider, for example, the case $c' = \dfrac{2-b_0}{d_0} < 1$ and
$< \dfrac{3-\mu_p}{m_p}$. In this case, we show that 
$p \notin \Bs|K_S+B+M-C_0-N|$, for some effective divisor $N$ supported away
from $p$. Using the vanishing 
\[
  H^1(S, K_S+\rup{M-c'D}) = H^1(S, K_S+B+M-2C_0-N) = 0,
\]
it suffices to show that $p \notin \Bs \left| K_S+B+M-C_0-N |_{C_0} \right|$;
this, in turn, will follow from (1.2.2)(a) and $\rup{M-c'D} \cdot C_0 > 1$.

Now $C_0$ passes through $p$ but is not tangent to $\vv$, and therefore we
have only $M \cdot C_0 \geq \beta_1$ (rather than $2\beta_1$). 
$\rup{M-c'D} = (M-c'D) + \Delta'$, with $\Delta' = B+c'D-2C_0$ in a
neighborhood of $p$, and therefore
$\ord_p(\Delta') \geq \frac{1}{2}(\mu_{\vv}+c'(4-\mu_{\vv})-2)$; 
it suffices to show that
\[
  (1-c't)\beta_1 + \frac{1}{2}(\mu_{\vv}+c'(4-\mu_{\vv})-2) > 1.
\]
An inequality equivalent to this one was already proved in (3.16).

The proof in the remaining cases is a similar adaptation of the arguments in
(3.13).
\end{proof}


\section{Example}

Fix an integer $n$, $n \geq 1$. Let $S$ be the $n^{\text{th}}$ Hirzebruch
surface, i.e. the geometrically ruled rational surface $\Prsp(\Ecal)$, where
$\Ecal$ is the rank 2 vector bundle $\Ocal_{\Prsp^1}\oplus\Ocal_{\Prsp^1}(-n)$
on $\Prsp^1$. Let $\pi:S\to\Prsp^1$ be the ruling of $S$, and let $F$ denote a
fiber of $\pi$. $S$ contains a unique irreducible curve $G$ with negative
self-intersection, $G^2=-n$. $\Pic(S) \cong \Z \oplus \Z$, with generators $F$
and $G$; $F^2=0, F \cdot G = 1$. $K_S \sim -2G-(n+2)F$. If $C$ is any
irreducible curve on $S$, then $C=G$, $C \sim F$, or $C \sim aG+bF$ with
$a,b \in \Z, a \geq 1$ and $b \geq na$. All these properties are proved, for
example, in \cite[Ch.V, \S 2]{harbook}.

\vspace{4pt}

Let $H_m=G+mF$. We will use the Reider-type results for \Q-divisors to prove
the following facts:

{\bf Claim.} (1) (See \cite[Ch.IV, Ex.1]{bv}.) 
\emph{
  $|H_n|$ is base-point-free, and defines a morphism $\phi_n:S\to\Prsp^{n+1}$.
  Moreover, $\phi_n$ is an isomorphism on $S \setminus G$, and
  $\bar{S} = \phi_n(S) \subset \Prsp^{n+1}$ is a (projective) cone with vertex
  $x = \phi_n(G)$. (Thus $\bar{S}$ is the cone over a normal rational curve
  contained in a hyperplane $\Prsp^n \subset \Prsp^{n+1}$, because
  $G \cong \Prsp^1$ and $G^2=-n$. $\phi_n$ is the blowing-up
  of $\bar{S}$ at $x$.)
}

(2) \emph{
  $|H_m|$ is very ample for $m \geq n+1$, defining an embedding 
  $\phi_m : S \to \Prsp^{2m-n+1}$.
}
(See \cite[Ch.IV, Ex.2]{bv} for other properties of $|H_m|$.)

\vspace{6pt}

Certainly, these facts can be proved in many different ways. For
instance, for (1): if $G \subset S$ is a smooth rational curve with negative
self-intersection on \emph{any} smooth surface $S$, then there is a projective
contraction $\phi : S \to \bar{S}$, which is an isomorphism on $S
\setminus G$ and contracts $G$ to a normal point $x$. (This is a direct
generalization of the ``easy'' part of Castelnuovo's criterion -- the ``hard''
part being the regularity of $\bar{S}$ at $x$ when $G^2=-1$.) The proof can be 
adapted to the situation of the Claim.  Alternatively,
most of the Claim is proved in \cite[Theorem 2.17]{harbook}.

The methods used in these proofs are somewhat specialized (the ``normal
contraction'' approach depends on $\Pic(G) \cong \Z$; the proof in
\cite{harbook} is typical for ruled surfaces). From this point of view,
Reider's theorem, which is based only on intersection numbers, is much more
general. However, as we will see, Reider's theorem doesn't apply in the
situation of the Claim. The proof we give below shows that there are instances
where the scope of Reider's original results can be broadened by allowing
\Q-divisors into the picture.

\newpage

\emph{Proof of the Claim.}

(1) Write $H_n=K_S+L$, thus defining $L = H_n-K_S = 3G+(2n+2)F$. Then
$L \cdot G = 3(-n)+(2n+2)=-n+2$; thus $L$ is not nef for $n \geq 3$, and
therefore Reider's criterion does not apply. However, write $L$ as $B+M$, with
$B=(1-\epsilon)G$ and $M=(2+\epsilon)G+(2n+2)F,\,\epsilon \in (0,1)$. Then
\[
  M \cdot F = (2+\epsilon), \quad M \cdot G = (2-\epsilon n), \quad
    M^2 = (2+\epsilon)(2n+4-\epsilon n).
\]
In particular, for $\epsilon \to 0$, we have $M \cdot F \to 2, M \cdot G \to 2, 
M^2 \to 2(2n+4) \geq 12$. Fix $\epsilon > 0, \epsilon \ll 1$, such that
$M^2>9, M \cdot F \geq \frac{3}{2}, \text{ and }M \cdot G \geq \frac{3}{2} $.
Since \emph{any} irreducible curve $C \subset S$ is either $C=G$, or $C \sim
F$, or $C \sim aG+bF$ with $a \geq 1$ and $b \geq na$, we automatically have 
$M \cdot C \geq \frac{3}{2}$ for all such $C$. (We will use this observation
again later: if $M\cdot G \geq 0$ and $C \neq G$ is an irreducible curve,
then $M \cdot C \geq M \cdot F$.) Therefore $|H_n|$ is base-point-free
by Proposition 3, part 2, with $\beta_2 = 3$ and $\beta_1 = \dfrac{3}{2} =
\dfrac{\beta_2}{\beta_2 - 1}$.

\vspace{4pt}

Thus $|H_n|$ defines a morphism $\phi_n:S\to\Prsp^{\nu}$, $\nu = \dim |H_n|$.
We compute $\nu$. By Riemann--Roch, we have:
\[
  \chi(S,H_n) = \frac{H_n \cdot (H_n-K_S)}{2} + \chi(S,\Ocal_S) = n+2.
\]
We get $\nu = h^0(S,H_n)-1=n+1$, as stated in the Claim, \emph{if} we can
show that $h^i(S,H_n)=0$ for $i\geq 1$. If we write $H_n=K_S+L$, as before,
Kodaira's vanishing theorem does not apply, because $L$ is not ample (it is
not even nef). If we write $L=B+M$ as above, though, we get $h^i(S,H_n)=0$
for $i \geq 1$, by (1.2.1)(a).

\vspace{4pt}

Next we show that $\phi_n$ is an isomorphism on $S \setminus G$. Consider
two distinct points $p,q \in F \setminus G$. Write $L=B'+M'$, with
$B'=(1-\epsilon)G+(1-\alpha)F, \,M'=(2+\epsilon)G+(2n+1+\alpha)F, \;
\epsilon, \alpha \in (0,1)$. (Note that we may use \emph{any} decomposition of
$L$ of the form $B+M$, as long as $\rup{M}=L$.) We have:
\begin{gather*}
   M' \cdot F = (2+\epsilon), \quad M' \cdot G = (1+\alpha-\epsilon n), \\
     (M')^2 = (2+\epsilon)(2n+2+2\alpha - \epsilon n).
\end{gather*}
In particular, for $\epsilon, \alpha \to 0$, $M'$ is nef and big and
$M' \to 2(2n+2) \geq 8$. Let $\mu \eqdef \mu_p = \mu_q = 1-\alpha$.
Choose $\beta_2 = \beta_{2,p} = \beta_{2,q} = \frac{3}{2}$ (say); then
$\beta_2 \geq 2-\mu = 1+\alpha$ and $(M')^2 > 2(\beta_2)^2$ for 
$\epsilon, \alpha \ll 1$. 

Fix $\epsilon \ll 1$, and then choose $\alpha \ll \epsilon$ such that 
$1+\dfrac{\epsilon}{2} \geq \dfrac{\beta_2}{\beta_2-(1-\mu)} =
\dfrac{\beta_2}{\beta_2 - \alpha}$; this can be done, because 
$\dfrac{\beta_2}{\beta_2-\alpha} \to 1$ for $\alpha \to 0$. Then
$M' \cdot F = 2+\epsilon = 2\beta_1$, with $\beta_1=1+\dfrac{\epsilon}{2}$ ---
and therefore $M' \cdot C \geq 2\beta_1$ for every irreducible curve $C$
through $p$ or $q$. Hence Proposition 4, part 3, applies (with $\beta_{1,p} =
\beta_{1,q} = \beta_1$): $|H_n|$ separates $(p,q)$.

If $p,q \in S \setminus G$ are distinct points on another irreducible
curve $\bar{F} \sim F$, the proof is similar --- take $B' = (1-\epsilon)G +
(1-\alpha)\bar{F}$. (We say $\bar{F} \sim F$ instead of ``fiber of $\pi : S 
\to \Prsp^1$'', to emphasize that the proof uses numerical arguments only.)
Finally, if no such curve passes through both $p$ and $q$, the proof is even 
easier.

Separation of tangent directions on $S \setminus G$ is proved exactly
the same way; note that $\mu_p(B') = \mu_V(B') = 1-\alpha$ if $B' =
(1-\epsilon)G+(1-\alpha)F,\, p \in F \setminus G, \text{ and } \vv \in
T_p(F) \setminus \{ \vec{0} \}$.

$H_n\cdot G=0\text{ and } H_n\cdot F=1$; therefore $\phi_n$ contracts $G$ to a
point $x \in \bar{S} = \phi_n(S) \subset \Prsp^{n+1}$, and $\phi_n(\bar{F})$
is a straight line in $\Prsp^{n+1}$ for every $\bar{F}\sim F$.

\vspace{8pt}

(2) As in part (1) of the Claim, we can show that $|H_m|$ is base-point-free
for $m \geq n+1$, and defines a morphism $\phi_m : S \to \Prsp^{2m-n+1}$ which
is an isomorphism on $S \setminus G$. For $m \geq n+1$, we must show that
$|H_m|$ separates $p,q$ even when $p$ (or $q$, or both) is on $G$, and also
that $|H_m|$ separates tangent directions at every point $p \in G$.

Let $\{p\}=F \inters G\text{ and }\vv\in T_p(G) \setminus\{ \vec{0} \}$. 
We will show that $|H_{n+1}|$ separates $\vv$ at $p$; the other properties
have similar proofs.

Write $H_{n+1} = K_S+L,\, L = 3G+(2n+3)F$. Write $L = B+M, \,
B=(1-\epsilon)G, \, M = (2+\epsilon)G+(2n+3)F, \; \epsilon \in (0,1)$.
We have:
\[
  M \cdot F = (2+\epsilon), \quad M \cdot G = (3-\epsilon n), \quad
    M^2 = (2+\epsilon)(2n+6- \epsilon n).
\]
For $\epsilon \to 0$ we have $M \cdot F \to 2, \, M \cdot G \to 3, \, 
\text{ and } M^2 \to 2(2n+6) \geq 16$; in particular $M$ is nef and big.
(Note that $L$ itself is not nef, if $n \geq 4$; indeed, $L \cdot G = 3-n$.)
We have $M \cdot C \geq 2+\epsilon$ for every irreducible curve $C \subset S$
(assuming that $\epsilon \ll 1$); also, if $\vv \in T_p(C)$, then 
$M \cdot C \geq 3 - \epsilon n$, because in that case $C \sim aG+bF$ with
$a \geq 1$ (proof: if $C \neq G$, then $C \cdot G \geq 2$, because
$\vv \in T_pC \inters T_pG$; therefore $C \not\sim F$.)

We have $\mu_p = \mu_V = 1-\epsilon$, and $\mu_{\vv} = 2(1-\epsilon)$. Choose
$\beta_2 = \beta_{2,p} = \beta_{2,V} = 2$ (say), so that $M^2 > 2(\beta_2)^2,
\, \beta_{2,p} \geq 2-\mu_p, \, \text{ and } \beta_{2,V} \geq 2-\mu_V$. 
Put $\beta_1 = \dfrac{2\beta_2}{2\beta_2-(2-\mu_{\vv})} = 
\dfrac{\beta_2}{\beta_2-\epsilon}$. For $\epsilon \ll 1$, we have: 
\begin{align*}
  M\cdot C &= 2+\epsilon\geq\beta_1 \quad\text{for all curves $C\subset S$}, \\
  M\cdot C &= 3-\epsilon n \geq 2\beta_1 \quad \text{ for all $C$ containing
  		$(p,\vv)$}.
\end{align*}
(Note that $\beta_1 = \dfrac{\beta_2}{\beta_2-1} \to 1$ as $\epsilon \to 0$, 
so these inequalities are verified for all $\epsilon \ll 1$.) Now use
Proposition 5, case 3.						\qed

\vspace{8pt}

By inspecting the proof of the Claim, we can see that the only assumptions
we used were that $\Pic(S) = \Z G \oplus \Z F, \, G^2=-n, \, F^2=0, \,
G \cdot F = 1, \, \text{ and } K_S=-2G-(n+2)F$ (if the other hypotheses are
satisfied, the last condition is equivalent to: $G$ and $F$ are smooth rational
curves); this suggests the following

\vspace{4pt}

{\bf Exercise.} A surface $S$ with these properties is isomorphic to
the $n^{\text{th}}$ Hirzebruch surface.

\vspace{4pt}

\emph{Hint.} There are several ways to see this. One, of course, is to use
part (1) of the Claim: after all, we have shown that $S$ is the blowing-up of
the cone over the normal rational curve of degree $n$.

Another solution is to show that $|F|$ is base-point-free and $\dim |F| = 1$,
as in the proof of part (1) of the Claim; thus $\phi = \phi_{|F|}$ realizes $S$
as a geometrically ruled surface over $\Prsp^1$, as required. ($S$ is minimal,
because $C^2 \geq 0$ for every irreducible curve $C \neq G$; this follows
easily from the hypotheses.)





\end{document}